\documentclass[reqno]{amsart}

\usepackage{amssymb,amsmath,enumerate}
\usepackage{graphics}
\usepackage{graphicx}
\usepackage{caption}
\usepackage{subcaption}
\usepackage[linesnumbered,boxruled,commentsnumbered]{algorithm2e}
\usepackage{comment,color}
\usepackage{hyperref}

\DeclareMathOperator*{\argmin}{arg\,min}

\newcommand{\Norm}[1]{\left|\left| #1 \right|\right|}

\def\R{\mathbb{R}}

\newcommand\bq{\mathbf{q}}
\newcommand\bp{\mathbf{p}}

\newcommand\bx{\mathbf{x}}

\newcommand\by{\mathbf{y}}

\newcommand\calX{\mathcal{X}}
\newcommand\calY{\mathcal{Y}}

\newcommand\thetab{\boldsymbol{\theta}}

\numberwithin{equation}{section}

\graphicspath{{./figures/}{.}}

\begin{document}

\title[Convex PINNS for Monge-Amp\`ere Equation]{Convex Physics Informed Neural Networks for the Monge-Amp\`ere Optimal Transport Problem}


\author[A. Caboussat]{Alexandre Caboussat}

\address{
Geneva School of Business Administration (HEG-Gen\`eve), 
University of Applied Sciences and Arts Western Switzerland (HES-SO), 
1227 Carouge, Switzerland \newline
Email : \texttt{alexandre.caboussat@hesge.ch}
}

\author[A. Peruso]{Anna Peruso}

\address{
Institute of Mathematics, Ecole Polytechnique F\'ed\'erale de Lausanne, 
1015 Lausanne, Switzerland, 
Email : \texttt{anna.peruso@epfl.ch}
and
Geneva School of Business Administration (HEG-Gen\`eve), 
University of Applied Sciences and Arts Western Switzerland (HES-SO), 
1227 Carouge, Switzerland, 
Email : \texttt{anna.peruso@hesge.ch} 
}


\date{\today}

\begin{abstract}
Optimal transportation of raw material from suppliers to customers is an issue arising in logistics that is addressed here with a continuous model relying on optimal transport theory.
A physics informed neural network method is advocated here for the solution of the corresponding generalized Monge-Amp\`ere equation. 
Convex neural networks are advocated to enforce the convexity of the solution to the Monge-Amp\`ere equation and obtain a suitable approximation of the optimal transport map. 
A particular focus is set on the enforcement of transport boundary conditions in the loss function. 
Numerical experiments illustrate the solution to the optimal transport problem in several configurations, and sensitivity analyses are performed. 
\end{abstract}


\maketitle
 
\section{Introduction}

In the field of fully nonlinear equations, the optimal transport problem has raised numerous fundamental questions over the recent years \cite{Ambrosioetal2003,PhilippisFigalli2014,Villani2003}. In particular, it has been linked with analytical models, such as, e.g., the Monge-Amp\`ere equation \cite{CaffarelliCabre1995,Figalli2017,LakkisPryer} or the prescribed Jacobian equation \cite{DacorognaMoser1990}. 
In the spirit of Monge \cite{Monge1781}, the optimal transport problem describes the problem of moving bulk material from a source (supplier) to a target destination location (customer) \cite{Bottanietal2019,oberman2015efficient,peyre2017computational}. 

\medskip

From the numerical viewpoint, many numerical methods from classical scientific computing have been advocated in the literature, involving, e.g., finite differences approaches \cite{Benamou2014NumericalSO,Benamouetal2010,Cassini2024} and finite element methods \cite{BenamouBrenier2000,PrinsSISC2015,Yadav2019}.
On the other hand, the emergence of deep neural networks has opened many avenues of investigation to approximate partial differential equations \cite{e2018deepritz,HORNIK1989359,VarNet2020}. 
In particular, Physics Informed Neural Networks (PINNs) have shown their capabilities to approximate solutions of partial differential equations \cite{RAISSI2019686}. 
Several extensions of PINNs are currently being investigated, e.g., consistent PINNs \cite{bonito2024} or variational PINNS \cite{Kharazmi2019VariationalPN}. 
Applications of PINNs to the classical Monge-Amp\`ere equation (with Dirichlet boundary conditions) can be found in \cite{NYSTROM}, while a first tentative to solve the optimal transport problem is presented in \cite{singh2022physics}. 

\medskip

We aim here at solving the following optimal transport problem. Let $\mathcal{X}$ and $\mathcal{Y}$ be two bounded, smooth open sets in $\mathbb{R}^2$, equipped respectively with probability measures $\mu(d\bx) = f(\bx)\,d\bx$ and $\nu(d\by) = g(\by)\,d\by$, where the density functions $f$ and $g$ are supported on $\mathcal{X}$ and $\mathcal{Y}$, respectively. We further assume that $f$ and $g$ are strictly positive and bounded on their respective supports. Brenier's Theorem \cite{Figalli2017}, implies that the optimal transport map, with respect to a quadratic cost, is given by the gradient of a convex potential $u$, which satisfies the generalized Monge-Ampère equation: find $u: \calX \rightarrow \mathbb{R}$ such that 
$$
\text{det}(D^2 u(\bx)) = \frac{f(\bx)}{g(\nabla u(\bx))} \quad \text{for } \mu\text{-almost every } \bx \in \mathcal{X},
$$
with the boundary condition
$$
\nabla u(\mathcal{X}) = \mathcal{Y}.
$$
To approximate the optimal transport map, it is important to ensure that the numerical method provides an approximation of $u$ that is also convex. 
Otherwise, there is no guarantee that the solution is unique. 
To guarantee that a neural network approximation of the solution remains convex, so-called {\em Input Convex Neural Networks} have been designed in the literature  \cite{amos17b}. 
First results for ICNNs applied to the optimal transport problem have been reported in \cite{singh2022physics}; they are specifically applied to the context of random sampling and density estimation based on available data, but do not aim at solving explicitly the partial differential equation. 
To ensure that the neural network approximation is convex, they require the activation function to be convex and non-decreasing, and all coefficients of the network to remain positive. 

\medskip

In this work, we propose to use ICNNs to approximate the solution of the generalized Monge-Ampère equation with transport boundary conditions. To enforce these boundary conditions at collocation points on the boundary, we incorporate a discretized version of the transport boundary condition into the loss function, inspired by the Hausdorff distance. 

The structure of this article is as follows: Section~\ref{sect:model} describes the derivation of the generalized Monge-Amp\`ere equation. Section~\ref{sect:pinns} describes the approach based on physics informed neural networks (PINNs), with, in particular ICNNs. Numerical experiments are reported in Section~\ref{sect:num} to validate the approach and highlight the role of input convex PINNs for this problem. 

\section{Mathematical Model}\label{sect:model}

\begin{figure}[ht!]
\begin{center}
\includegraphics[width=0.8\textwidth]{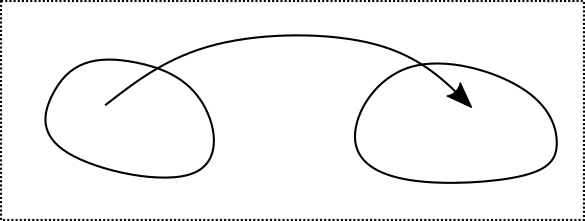}
\end{center}
\setlength{\unitlength}{1cm}
\begin{picture}(0,0)
\put(-4.8,4){$\Omega$}
\put(-3.8,2.2){$x$}
\put(-2.2,2){$\mathcal{X}$}
\put(-4,1.3){$\partial \mathcal{X}$}
\put(3,2){$\mathcal{Y}$}
\put(3,1){$\partial \mathcal{Y}$}
\put(-1,4){$\mathbf{m}(x)$}
\end{picture}
\caption{Sketch of the optimal transport problem, illustrating the mapping of a domain $\mathcal{X}$ onto another domain $\mathcal{Y}$.}
\label{fig1}
\end{figure}

The optimal transport problem, initially introduced by Monge \cite{Monge1781}, is formulated as follows: given two probability measures $\mu$ and $\nu$ defined on the measurable spaces $\mathcal{X}$ and $\mathcal{Y}$, find a measurable map $\mathbf{m}: \mathcal{X} \rightarrow \mathcal{Y}$ such that $\mathbf{m}_{\#} \mu = \nu$, i.e.,
\begin{equation}
\nu(A) = \mu \left( \mathbf{m}^{-1}(A) \right) \quad \forall \, A \subset \mathcal{Y} \text{ measurable},
\label{eq:pushforward}
\end{equation}
in such a way that $\mathbf{m}$ minimizes the transportation cost:
\begin{equation}
\int_{\mathcal{X}} c(\bx, \mathbf{m}(\bx)) \, d\mu(\bx) = \min_{\mathbf{s}_{\#} \mu = \nu} \int_{\mathcal{X}} c(\bx, \mathbf{s}(\bx)) \, d\mu(\bx),
\label{eq:cost_minimization}
\end{equation}
where $c: \mathcal{X} \times \mathcal{Y} \rightarrow \mathbb{R}$ is the quadratic cost function $c(\bx,\by) = \frac{1}{2}|\bx-\by|^2$, and the minimum is taken over all measurable maps $\mathbf{s}: \mathcal{X} \rightarrow \mathcal{Y}$ such that $\mathbf{s}_{\#} \mu = \nu$. When the condition $\mathbf{m}_{\#} \mu = \nu$ is met, $\mathbf{m}$ is called a transport map, and if $\mathbf{m}$ also minimizes the cost, it is referred to as an optimal transport map. 
Figure~\ref{fig1} provides a visual sketch of this scenario.

Even for the quadratic cost, determining the existence of an optimal transport map is not straightforward. Additionally, examples can be constructed where the Monge problem is ill-posed due to the non-existence of a transport map; for instance, when $\mu$ is a Dirac measure and $\nu$ is not. This suggests that certain conditions on the measures $\mu$ and $\nu$ are necessary. Consider the case where $\mathcal{X}, \mathcal{Y} \subset \mathbb{R}^n$, with $\mu(d\bx) = f(\bx) \, d\bx$ and $\nu(d\by) = g(\by) \, d\by$. If $\mathbf{s}: \mathcal{X} \rightarrow \mathcal{Y}$ is a sufficiently smooth transport map, the condition $\mathbf{s}_{\#} \mu = \nu$ can be reformulated as a Jacobian equation. Specifically, if $\chi: \mathbb{R}^n \rightarrow \mathbb{R}$ is a test function, the condition in \eqref{eq:pushforward} leads to
\[
\int_{\mathbb{R}^n} \chi(\mathbf{s}(\bx)) f(\bx) \, d\bx = \int_{\mathbb{R}^n} \chi(\by) g(\by) \, d\by.
\]
Assuming $\mathbf{s}$ is a diffeomorphism, we can set $y = \mathbf{s}(x)$ and apply a change of variables that yields
\[
\int_{\mathbb{R}^n} \chi(\mathbf{s}(\bx)) f(\bx) \, d\bx = \int_{\mathbb{R}^n} \chi(\mathbf{s}(\bx)) g(\mathbf{s}(\bx)) \det(\nabla \mathbf{s}(\bx)) \, d\bx.
\]
By the arbitrariness of $\chi$, this yields the Jacobian equation: $\mathbf{s}: \calX \rightarrow \mathbb{R}^n$ such that  
\[
f(x) = g(\mathbf{s}(\bx)) \det(\nabla \mathbf{s}(\bx)) \quad \text{a.e.} \ \bx \in \calX.
\]
In the specific case where $\mu$ and $\nu$ are absolutely continuous probability measures with compact support in $\mathcal{X}$ and $\mathcal{Y}$, respectively, with density functions $f$ and $g$ that are bounded away from zero and infinity on $\mathcal{X}$ and $\mathcal{Y}$, respectively, Brenier's theorem states that:

\begin{enumerate}[i)]
    \item There exists a unique solution $\mathbf{m}$ to the optimal transport problem \eqref{eq:cost_minimization}.
    \item There exists a convex function $u:\calX \rightarrow \mathbb{R}$ such that the optimal map $\mathbf{m}$ is given by $\mathbf{m}(\bx) = \nabla u(\bx)$ for $\mu$-almost every $\bx$, and satisfies
\begin{equation*}
\label{eq:Jacobian}
\text{det}(D^2u(\bx)) = \frac{f(\bx)}{g(\nabla u(\bx))} \quad \text{for } \mu\text{-a.e. } \bx \in \mathbb{R}^n.
\end{equation*}
\end{enumerate}
Moreover, the condition in \eqref{eq:pushforward} implies that the boundary condition $\nabla u(\mathcal{X}) = \mathcal{Y}$ 
can be replaced by the boundary condition
\begin{equation}\label{bc1}
\nabla u(\partial \mathcal{X}) = \partial \mathcal{Y}.
\end{equation}
In the sequel, let us now address the approximation of the optimal transport problem with neural networks.

\section{Physics Informed Neural Networks}\label{sect:pinns}

\subsection{Generalities}

We consider feedforward neural networks made up of an input layer, an output layer and $L\geq 1$ hidden layers. We denote by $N_j$ the number of neurons of the $j^{th}$ layer, $j=0,\ldots,L+1$. 
The network is parametrised by the weight matrices $\mathbf{W}^{(l)}\in\mathbb{R}^{N_{l+1}\times N_{l}}$ and bias vectors $\mathbf{b}^{(l)}\in\mathbb{R}^{N_{l+1}}$ with $l=0, \ldots, L$. 
The output is defined with respect to the input as: 
    \begin{align*}
    &\mathbf{x}^0 = \mathbf{x},\\
    &\mathbf{x}^l = \sigma(\mathbf{W}^{(l-1)}\mathbf{x}^{l-1}+\mathbf{b}^{(l-1)}), \quad 1\leq l\leq L\\
    &\mathbf{y} = \mathbf{W}^{(L)}\mathbf{x}^L+\mathbf{b}^{(L)},
    \end{align*}
where $\sigma$ is the activation function. 
The set of parameters is collectively denoted by $\thetab =\big\{\mathbf{W}^{(l)},\mathbf{b}^{(l)}\big\}_l$, and the neural network mapping is denoted by $\mathbf{y} := u_{NN} (\mathbf{x})$. 
Common choices for $\sigma$ include the ReLU function, $\sigma(z) = \max(0, z)$, the hyperbolic tangent function, $\sigma(z) = \frac{e^{z}-e^{-z}}{e^{z}+e^{-z}}$, or the softplus function $\sigma(z) = \log(1+e^z)$. 

A neural network has the potential to provide an approximation $u_{NN}:\mathbb{R}^{N_0}\to\mathbb{R}^{N_{L+1}}$ to a given mapping $u:\mathbb{R}^{N_0}\to\mathbb{R}^{N_{L+1}}$. 
Indeed, the universal approximation theorem \cite{HORNIK1991251} shows that, for single-layered neural networks and regardless of the shape of $\sigma$, this approximation is possible by density of the approximation spaces when the number $n$ of hidden neurons increases. 
Universal approximation theorems also apply to deep neural networks \cite{HORNIK1991251,HORNIK1989359} under some smoothness conditions.

\subsection{Input Convex Neural Networks}
As discussed in Section~\ref{sect:model}, the solution $u$ to the generalized Monge-Amp\`ere equation must be a convex function. 
Therefore, it is imperative to ensure the convexity of the approximate solution $u_{NN}$. 
Neural networks that enforce convexity are mandatory, and may come at a higher (computational) price. 
To this end, one can employ \emph{Input Convex Neural Networks} (ICNNs), as introduced in \cite{amos17b}, whose architecture inherently guarantees convexity with respect to the input. 
The general structure of an ICNN with $L$ hidden layers and input variables $x$ is given by
\begin{equation}
    \begin{aligned}
    \mathbf{x}^0&= \mathbf{x}, \\
    \mathbf{x}^1&= \sigma (\mathbf{L}^{(0)}\mathbf{x}^0 + \mathbf{b}^{(0)}), \\
    \mathbf{x}^{l} &= \sigma (\mathbf{W}^{(l-1)}\mathbf{x}^{l-1} + \mathbf{L}^{(l-1)}\mathbf{x}^0 + \mathbf{b}^{(l-1)}), \quad 2\leq l\leq L , \\
    \mathbf{y}&=\mathbf{W}^{(L)}\mathbf{x}^{L} + \mathbf{L}^{(L)}\mathbf{x}^0 + \mathbf{b}^{(L)}.
    \end{aligned}
\end{equation}
Here, $\sigma$ denotes a specified increasing convex function, $\mathbf{W}^{(l)}\in\mathbb{R}^{N_{l+1}\times N_{l}}$ are matrices with nonnegative entries, $\mathbf{L}^{(l)}\in\mathbb{R}^{N_{l+1}\times N_0}$ are matrices, and $\mathbf{b}^{(l)}\in\mathbb{R}^{N_{l+1}}$ are vectors. The convexity of the network follows from the following principles:
\begin{enumerate}[i)]
    \item A linear combination of convex functions with positive coefficients remains convex.
    \item The composition $h \circ g$ of a convex function $g : \mathbb{R}^n \rightarrow \mathbb{R}$ with an increasing convex function $h : \mathbb{R} \rightarrow \mathbb{R}$ is also convex.
\end{enumerate}

A notable feature of ICNNs is the inclusion of \textit{passthrough} layers, which directly connect the input $\mathbf{x}$ to the hidden units in deeper layers. 
In traditional feedforward networks, such layers are unnecessary, as previous hidden units can be mapped to subsequent hidden units via an identity mapping. 
However, in ICNNs, the non-negativity constraint on the weights $\mathbf{W}^{(i)}$ restricts the use of hidden units that replicate the identity mapping. 
Thus, this additional passthrough is explicitly incorporated to overcome this limitation and is crucial for the representational capacity of the network \cite{Chen2018}.

In \cite{NYSTROM}, which addresses the approximation of the Dirichlet Monge–Ampère problem, the authors have tested the ICNN architecture up to dimension $4$. 
They have observed that, in theory, the method is applicable in any dimension. 
However, in practice, training scales poorly with the space dimension, which poses a significant challenge for optimal transport problems in high dimensions. 
The authors have also explored an alternative approach to enforcing convexity, based on penalizing the diagonal elements of the Hessian matrix. 
This method is viable in two dimensions but it does not impose convexity in a structural way. 
As a result, convergence to a convex solution may only occur asymptotically.

Therefore, following the conclusions in \cite{NYSTROM}, ICNNs remain the best networks to enforce the convexity of the solution, despite being hard to train. 

\subsection{Physics Informed Neural Networks}
To solve Monge-Amp\`ere equations using neural networks approximation, we use Physics-Informed Neural Networks (PINNs) based on the strong formulation of the partial differential equation \cite{RAISSI2019686}. 
Consider the Monge-Amp\`ere problem with Dirichlet boundary conditions: find $u: \mathcal{X} \rightarrow \mathbb{R}$ such that:
\begin{equation}\label{gen_ma_bc}
\left\{
\begin{array}{ll}
\det{{D}^2 u}(\bx) =  \dfrac{f(\bx)}{g (\nabla u(\bx))} \quad& \bx \in \mathcal{X},  \\
u(\bx) = h(\bx)  & \bx \in \partial \mathcal{X},
\end{array} 
\right.
\end{equation}
where $h$ is a prescribed function. The PINNs framework seeks a neural network approximation $u_{NN}$ of the solution $u$ to \eqref{gen_ma_bc} by minimizing a weighted sum of the residual $\mathcal{L}_{PDE}$ and the boundary loss $\mathcal{L}_b$, specifically:
$$
u_{NN} := \argmin_{\thetab}\mathcal{L}(\thetab) = \argmin_{\thetab} \left\{\mathcal{L}_{PDE}(\thetab) + C \mathcal{L}_b(\thetab)\right\},
$$
where 
\begin{align*}
    \mathcal{L}_{PDE} &= \left\|\det{D^2 u_{NN}} - \displaystyle \frac{f}{g (\nabla u_{NN})} \right\|_{L^2(\mathcal{X})}, \\
    \mathcal{L}_b &= \| u_{NN} - h\|_{L^2(\partial \mathcal{X})},
\end{align*}
and $C>0$ controls the weight of the boundary conditions. 
To numerically approximate $\mathcal{L}_{PDE}$ and $\mathcal{L}_b$, we select a set of collocation points $\{{\bp}_i\}_{i=1}^{N_c} \subset \mathcal{X}$ and boundary points $\{{\bq}_i\}_{i=1}^{N_b} \subset \partial \mathcal{X}$, and approximate the $L^2$ norms as follows:
\begin{align*}
    \mathcal{L}_{PDE} \approx& E_{PDE} := \frac{1}{N_c}\sum_{i=1}^{N_c} \left\vert \det{{D}^2 u_{NN}({\bp}_i)} - \frac{f({\bp}_i)}{g (\nabla u_{NN}({\bp}_i))} \right\vert^2,\\
    \mathcal{L}_b \approx &E_b := \frac{1}{N_b}\sum_{i=1}^{N_b}\vert u_{NN}({\bq}_i) - h({\bq}_i)\vert^2.
\end{align*}

When the Dirichlet boundary conditions in \eqref{gen_ma_bc} are replaced by the transport boundary conditions \eqref{bc1}, we consider the problem
\begin{equation}\label{gen_ma_ot}
\left\{
\begin{array}{ll}
\det{{D}^2 u}(\bx) =  \dfrac{f(\bx)}{g (\nabla u(\bx))} \quad& \bx \in \mathcal{X},  \\
\nabla u(\partial \calX) = \partial\calY.
\end{array} 
\right.
\end{equation}
The residual $E_b$ must be adjusted accordingly. 
To do so, we notice that 
$$
\nabla u(\partial \calX) = \partial\calY \iff d_H(\nabla u(\partial\calX), \partial \calY)=0,
$$
where $d_H$ is the Hausdorff distance with respect to the Euclidean metric, which is defined by 
$$
d_H(A, B) := \max \left\{ \sup_{\mathbf{a} \in A}  \mathrm{dist}(\mathbf{a}, B),\; \sup_{\mathbf{b} \in B}  \mathrm{dist}(A, \mathbf{b}) \right\},\quad A,B\subset\mathbb{R}^n,
$$
with $\mathrm{dist}(\mathbf{a}, B):=\min_{\mathbf{b}\in B}\|\mathbf{a}-\mathbf{b}\|_2$. 
The Hausdorff distance not only ensures that a point in $\partial \calX$ is mapped through $\nabla u$ to a point in $\partial \calY$ but also guarantees the surjectivity of the map $\nabla u$. 
On a discrete level we approximate $d_H$ in the following way: suppose we have two sets of boundary points 
$$
\{{\bx_i}\}_{i=1}^{N_{bx}} \subset \partial \mathcal{X},\quad \{{\by_i}\}_{i=1}^{N_{by}} \subset \partial \mathcal{Y},
$$  
the boundary loss $E_{OT}$ is then defined as
\begin{equation}
    E_{OT} = \frac{1}{N_{bx}}\sum_{i=1}^{N_{bx}}\mathrm{dist}\left(\nabla u ({\bx_i}), \{{\by_j}\}_{j=1}^{N_{by}}\right)^2 + \frac{1}{N_{by}}\sum_{i=1}^{N_{by}}\mathrm{dist}\left(\nabla u \big(\{{\bx_j}\}_{j=1}^{N_{bx}}\big), {\by_i}\right)^2,
    \label{eq:tildeEb}
\end{equation}
where the first term in \eqref{eq:tildeEb} controls the injectivity of the map $\nabla u_{NN}$ while the second term its surjectivity. Finally, we look for $u_{NN}$ convex such that
$$u_{NN} = \argmin_{\thetab}\left\{E_{PDE} + CE_{OT}\right\}.$$
 
The Hausdorff distance used in the loss function  to enforce the transport boundary condition can be computationally expensive. 
However, the method is general and may be applicable to a variety of situations. 
Furthermore, in some particular cases, the transport boundary condition admits a simpler formulation, for instance when mapping a square domain onto itself.
It can be shown that, in such a case, the optimal transport condition is equivalent to a Neumann condition. 
Therefore, the boundary terms of the loss function can be simplified and the PINNs solver can be significantly optimized.

\section{Numerical experiments}\label{sect:num}

Several numerical experiments are presented in order to validate the convergence and accuracy properties of the convex PINNs methodology. 

\subsection{Implementation details}
For the implementation of the PINN solver, the code has been developed in \texttt{pytorch} and it is available on request. 
The input convex neural networks used in the test cases consist of 4 hidden layers with 10 neurons per layer. 
The activation function employed is the softplus function, defined as $\sigma(x) = \log(1 + e^x)$, and is chosen because it is an increasing strictly convex function. 
This choice follows the recommendations provided in \cite{NYSTROM,singh2022physics}. 
To guarantee the positivity of the weight matrices $\mathbf{W}^{(l)} \in \mathbb{R}^{N_{l+1} \times N_{l}}$, the weights are initialized by squaring them element-wise, i.e., $w_{ij} \leftarrow w_{ij}^2$, as suggested in \cite{NYSTROM}. 
Note that the alternative approach of clipping the weights, such that $w_{ij} \geq 0$ for all $i,j$, following \cite{amos17b}, does not provide satisfactory results in our experiments. 

Following standard practice in PINNs \cite{pinnstraining}, we perform the initial training with the Adam optimizer before switching to the L‑BFGS optimizer to achieve more stable and accurate convergence.
Moreover, ICNNs are notoriously difficult to train \cite{bunne2022,richterpowell2021} and a meaningful initialization is required. In particular, the network parameters are initialized by training the network to approximate the identity map for the gradient of the solution, $\nabla u \approx \text{Id}$, inspired by the work in \cite{singh2022physics}. Training details are summarized in Table~\ref{tab:lbfgs_only}.
We have not experienced any difficulties due to the squaring of the coefficients in the initial optimization for any experiments presented hereafter. 
By default, every test case is run with $N_c = 800$ collocation points on the interior of the domain and $N_b = 800$ points on the boundary, choice which will be justified in the numerical experiments in the sequel. 
The constant $C$ is set to 1, unless specified otherwise. 

\begin{table}[htbp]
  \centering
  \caption{L‑BFGS training configuration.}
  \label{tab:lbfgs_only}
  \begin{tabular}{l|l}
    Main optimizer                  & L‑BFGS                                        \\
    \hline
    Line search                     & strong Wolfe                                 \\
    Sub-iterations per epoch        & 20                                           \\
    Learning rate                   & 1                                            \\
    Total epochs                    & 100                                          \\
    Batch size                      & Full-batch                                   \\
  \end{tabular}
\end{table}

\subsection{Transport of a constant density from a unit disk domain to an ellipse domain}
\begin{figure}[ht!]
\begin{center}    
\includegraphics[width=0.75\linewidth]{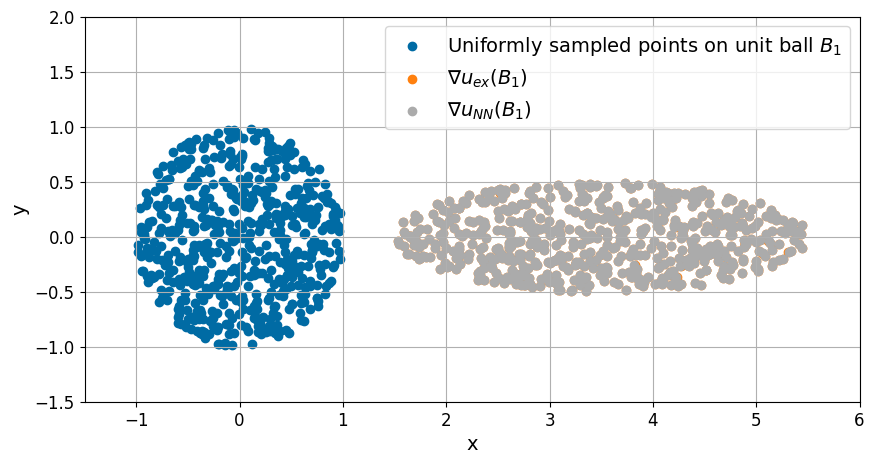} 
\end{center}    
    \caption{Transport of a constant density from a unit disk domain to an ellipse domain. 
    Visualization of uniformly sampled source points and their image through exact and approximate transport maps.}
    \label{fig:otfirst}
\end{figure}

\begin{figure}[ht!]
\begin{center}    
\includegraphics[width=0.5\linewidth]{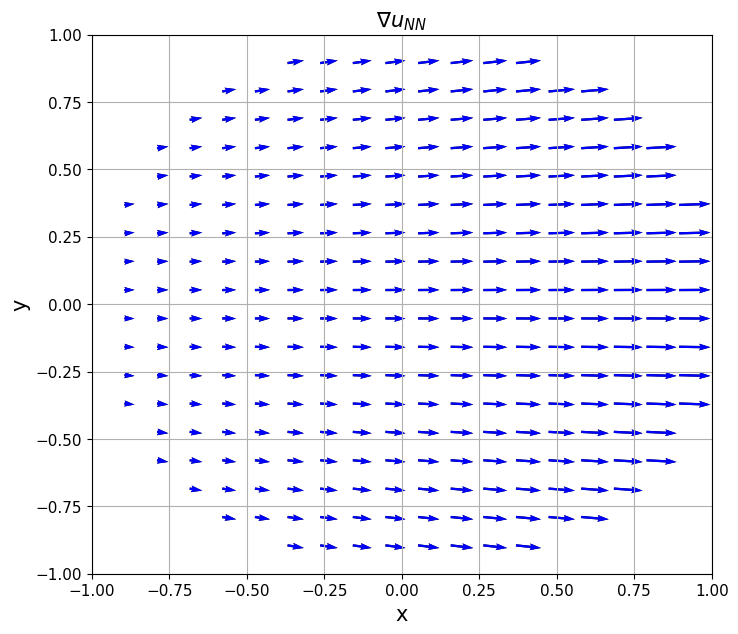} 
\end{center}    
\caption{Transport of a constant density from a unit disk domain to an ellipse domain. 
Visualization of the transport map $\nabla u_{NN}$.
}
    \label{fig:otfirst_field}
\end{figure}

\begin{figure}[ht!]
\begin{center}
\begin{tabular}{c}
\includegraphics[width=0.75\linewidth]{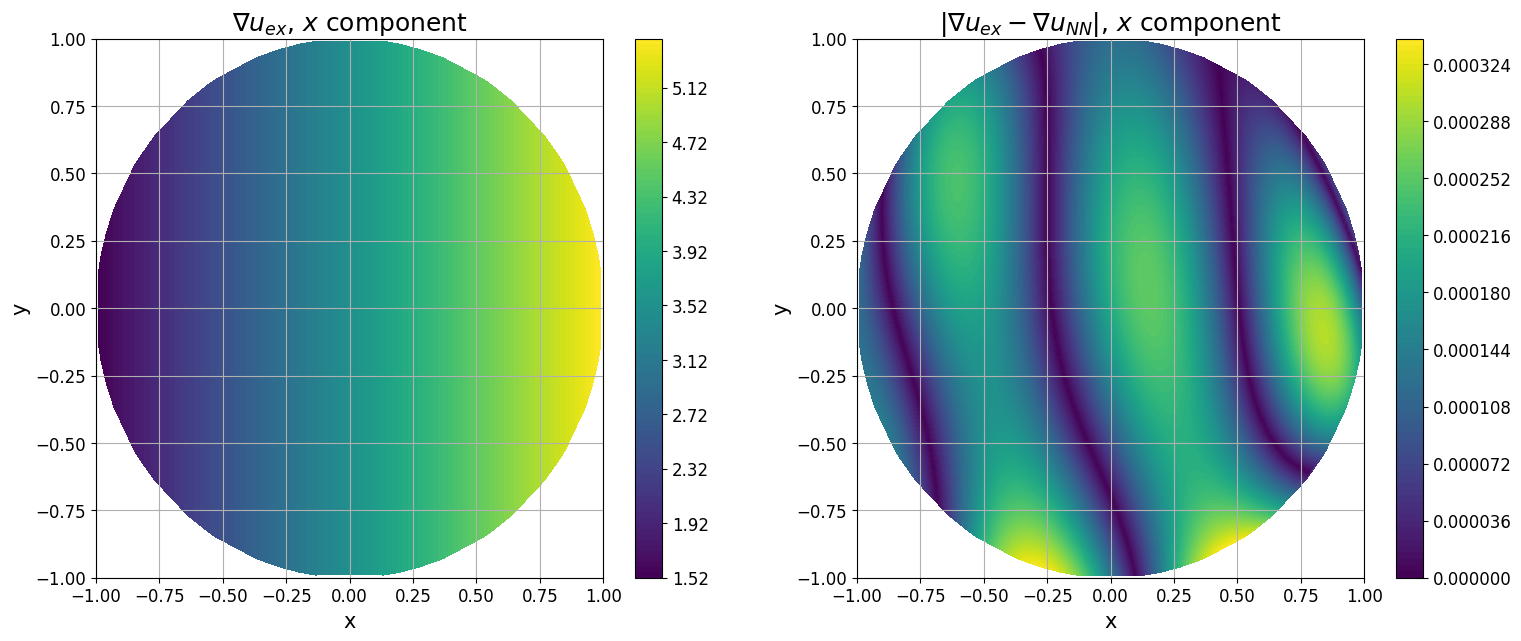} \\
\includegraphics[width=0.75\linewidth]{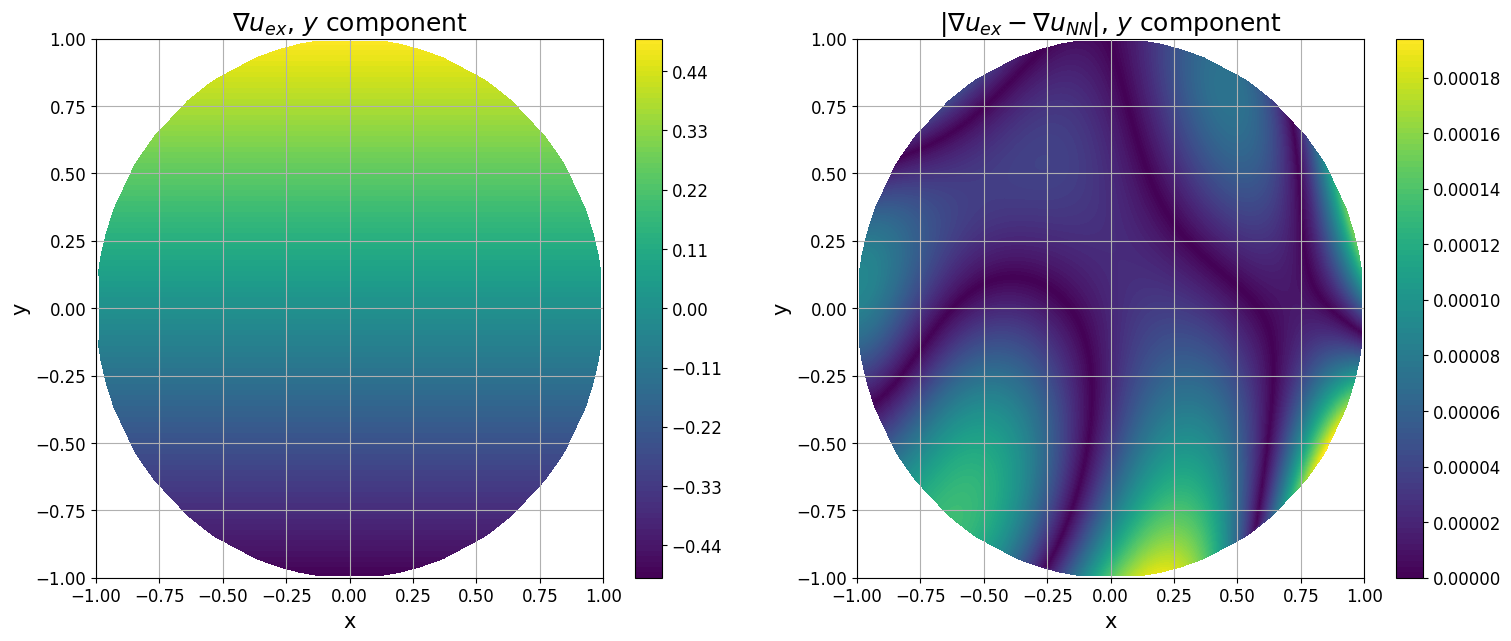} 
\end{tabular}
\end{center}    
    \caption{Transport of a constant density from a unit disk domain to an ellipse domain. 
    Top left: first $x$-coordinate of the transport map $\nabla u_{NN}$; 
    Top right: Pointwise absolute error $\left| \nabla u_{NN,x} - \nabla u_{ex,x} \right|$;  
    Bottom left: second $y$-coordinate of the transport map $\nabla u_{NN}$; 
    Top right: Pointwise absolute error $\left| \nabla u_{NN,y} - \nabla u_{ex,y} \right|$. }
    \label{fig:otfirst_xerr_yerr}
\end{figure}

\begin{figure}[ht!]
\begin{center}
\begin{tabular}{cc}
\includegraphics[width=0.45\linewidth]{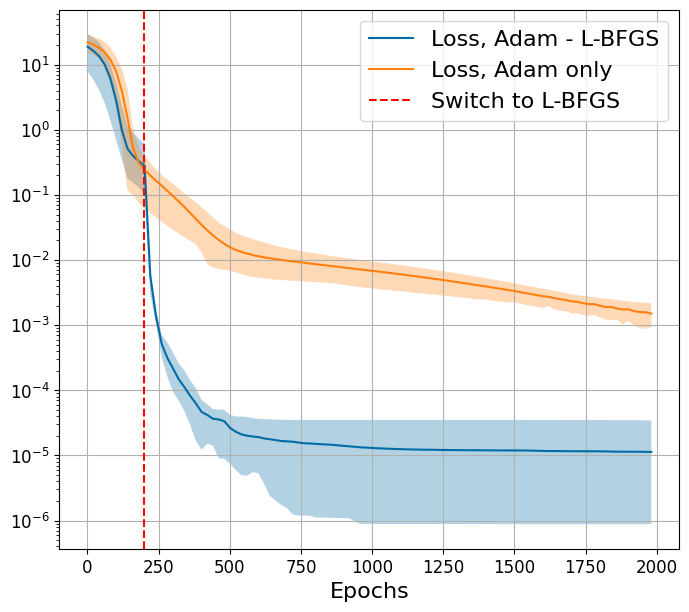}
\includegraphics[width=0.45\linewidth]{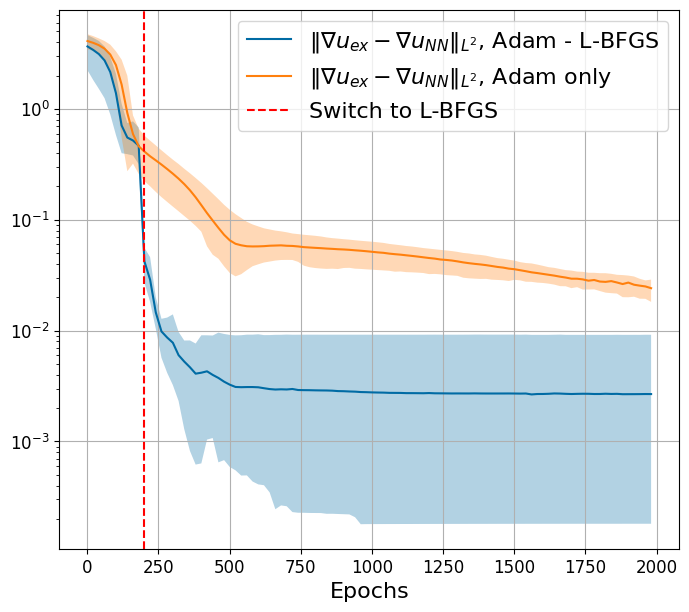} 
\end{tabular}
\end{center}    
    \caption{
    Transport of a constant density from a unit disk domain to an ellipse domain. 
    Left: Loss function vs epochs. 
    Right: $\|\nabla u - \nabla u_{NN}\|_{L^2}$ vs epochs, with the error computed on a validation set. 
    For the L-BFGS optimizer, the total epoch count is taken as the number of (outer) L‑BFGS iterations times the number of inner sub‐iterations. 
    Shaded bands show the 5th–95th percentile range over ten independent runs.
    }
    \label{fig:otfirst_loss}
\end{figure}

The first example aims at transporting a constant probability density function supported in a unit disk domain $ \calX \subset \R^2$ to another constant probability density function whose support is an ellipse domain $ \calY \subset \R^2$.
Let us thus define
$$ 
\mathcal{X}=\{\bx=(x_1,x_2)\in\R^2:\,  x_1^2+x_2^2 < 1\},
$$
and  
$$ 
\mathcal{Y} = \left\{\by=(y_1,y_2) \in \R^2:\, \frac{\left(y_1 - 3.5 \right)^2}{(2)^2}+\frac{y_2^2}{(0.5)^2} < 1\right\}.
$$ 
The corresponding generalized Monge-Amp\`ere problem (\ref{gen_ma_ot}) is considered here with 
$f(\bx) = \frac{1}{\pi}\chi_{\mathcal{X}}(\bx)$, $g(\by) = \frac{1}{\pi}\chi_{\mathcal{Y}}(\by)$, where $\chi_{A}$ denotes the characteristic function of the domain $A$. 
With this data, the exact solution to (\ref{gen_ma_ot}) is 

$$
\nabla u_{ex}(\bx) 
=\begin{pmatrix} 2x_1 + \frac{7}{2} \\ \frac{1}{2}x_2 \end{pmatrix}, 
\quad 
\bx = (x_1,x_2)\in \mathcal{X}.
$$
Figure~\ref{fig:otfirst} shows how uniformly sampled points in the unit ball are transported to the ellipse, both with the exact solution and the approximate map. 
It illustrates an appropriate, uniform, repartition of the sampled points (left), and an accurate transport of those points, whose final position coincides well with their exact values (right). 
Figure~\ref{fig:otfirst_field} visualizes the corresponding vector field $\nabla u_{NN}$, while Figure~\ref{fig:otfirst_xerr_yerr} details the first and second components of $\nabla u_{NN}$, together with the numerical error component-wise. 
We compare how the L-BFGS performs against the Adam optimizer. In order to do so, we run ten independent trainings, each of them with a different random draw of collocation and boundary points and a different random initialization of network parameters. 
Figure~\ref{fig:otfirst_loss} 
demonstrates that the L‑BFGS solver allows the loss function to decrease more rapidly but also to reach a lower final error than the Adam solver. 
Table~\ref{tab:lbfgs_vs_adam} 
reports the corresponding averages over these ten runs.

\begin{table}[htbp]
  \centering
  \caption{Comparison of the Adam vs. L‑BFGS optimizers. Provided values are the average over ten independent runs, (together with standard deviation).}
  \label{tab:lbfgs_vs_adam}
  \begin{tabular}{c|ccc}
    \textbf{Optimizer} & \textbf{Final Loss} & \boldmath$L^2$ \textbf{Error (Train)} & \boldmath$L^2$ \textbf{Error (Test)} \\
    \hline
    Adam   &  
		$0.00150\,(0.0043)$ &  
		$0.02399\,(0.00327)$ &
      	$0.02504\,(0.00389)$ \\
    L‑BFGS &  
      	$1.13\cdot 10^{-5}\,(1.39\cdot 10^{-5})$ &  
      	$0.00267\,(0.00376)$ &  
      	$0.00268\,(0.00376)$ \\
  \end{tabular}
\end{table}

\begin{figure}[ht!]
\begin{center}    
\begin{tabular}{c}
\includegraphics[width=0.58\linewidth]{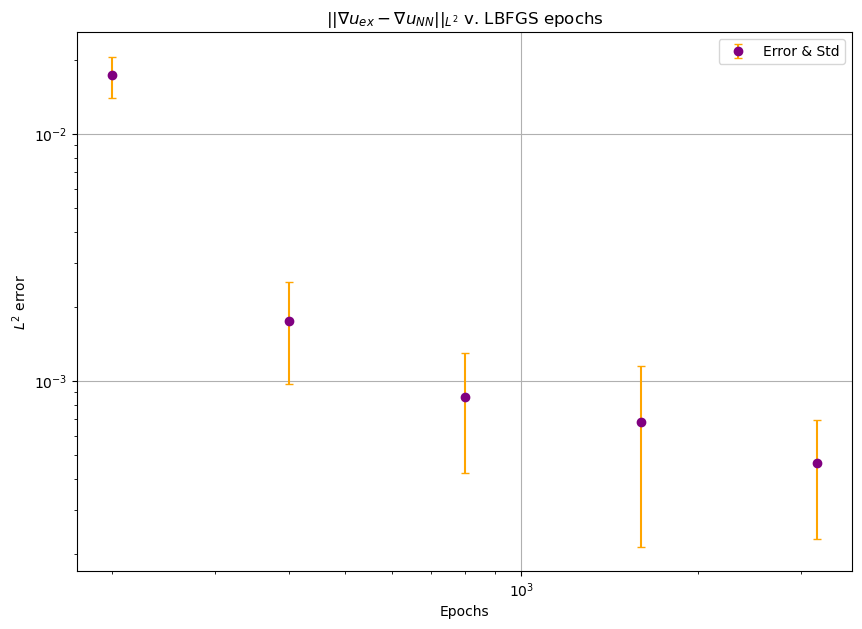} 
\\
\includegraphics[width=0.58\linewidth]{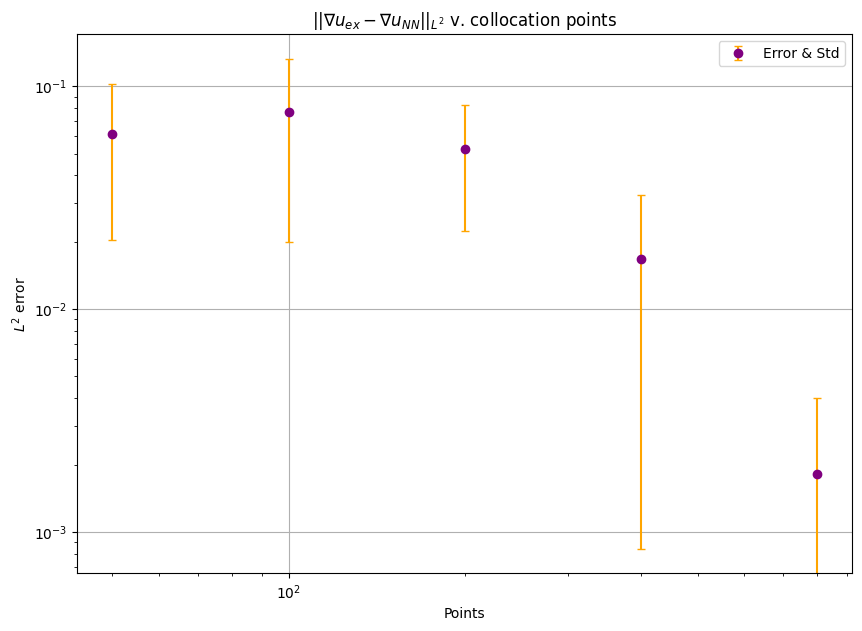} 
\\
\includegraphics[width=0.58\linewidth]{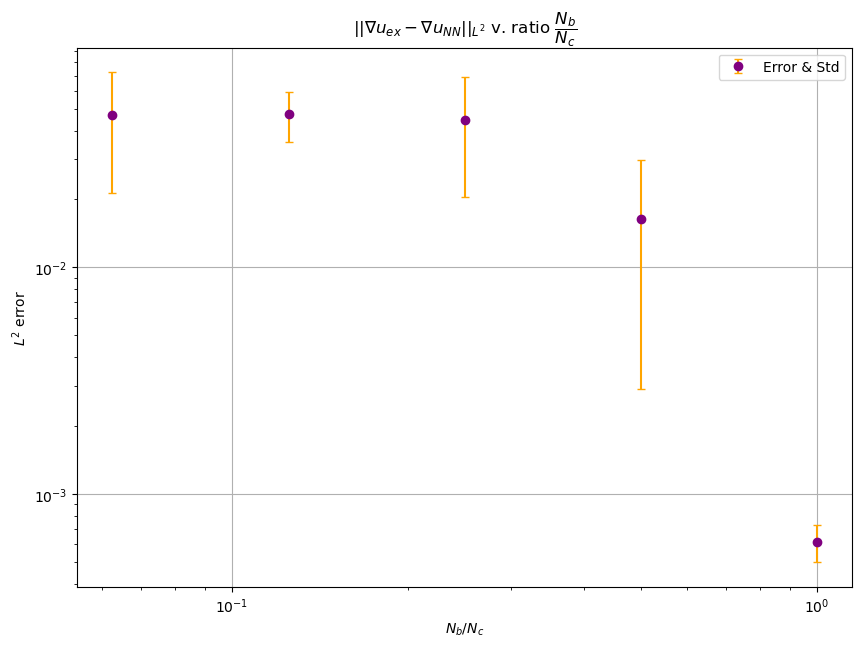} 
\end{tabular}
\end{center}
\caption{
Transport of a constant density from a unit disk domain to an ellipse domain. 
Top: $L^2$ (test) error $\Norm{\nabla u_{ex} - \nabla u_{NN}}_{L^2(\calX)}$ vs number of epochs. 
Middle: $L^2$ (test) error $\Norm{\nabla u_{ex} - \nabla u_{NN}}_{L^2(\calX)}$ vs number of collocation points. 
Bottom: $L^2$ (test) error $\Norm{\nabla u_{ex} - \nabla u_{NN}}_{L^2(\calX)}$ vs ratio $N_b/N_c$.
(Results averaged on $10$ simulations.)}
\label{fig:otfirst_sens}
\end{figure}

Let us study the sensitivity of those results with respect to the number of epochs and number of collocation points. 
Figure~\ref{fig:otfirst_sens} (top and middle) visualizes the behavior of the $L^2$ error $\Norm{\nabla u_{ex} - \nabla u_{NN}}_{L^2(\calX)}$ with respect to these factors respectively. 
As before, the error values specified are the values obtained when averaged over ten independent runs, each of them with a different random draw of collocation and boundary points and a different random initialization of network weights. 
In particular, when the number of collocation points increases, e.g., from $N_c=50$ to $N_c=100$, the first $50$ points in each sample remains the same as those for the smaller sample. 
It shows that the error decreases when either of these quantities increases, albeit with a large sample variation. 
Finally, Figure~\ref{fig:otfirst_sens} (bottom) visualizes the behavior of the $L^2$ error $\Norm{\nabla u_{ex} - \nabla u_{NN}}_{L^2(\calX)}$ with respect to the ratio $N_b/N_c$ between the number of boundary points and the number of collocation points. It shows that the error is smaller for a ratio equal to one, which confirms the initial choice of $N_b = N_c$ made earlier. 
More precisely, Figure~\ref{fig:otfirst_sens} (bottom) shows that the ratio between the number of collocation points, and the number of boundary points can be considered to be smaller or equal to one. 
However, considering a ratio below one corresponds to having more points on the boundary than in the interior of the domain. 
It is generally not meaningful to place more collocation points on the boundary than in the interior of the domain, since the boundary conditions represent lower-dimensional constraints (1D) compared to the governing equation defined in the 2D interior.
Thus, we will consider $N_b = N_c$ in the sequel.

\subsection{
Transport of a constant density from an ellipse domain to a rotated ellipse domain 
}

\begin{figure}[ht!]
\begin{center}
\includegraphics[width=0.6\linewidth]{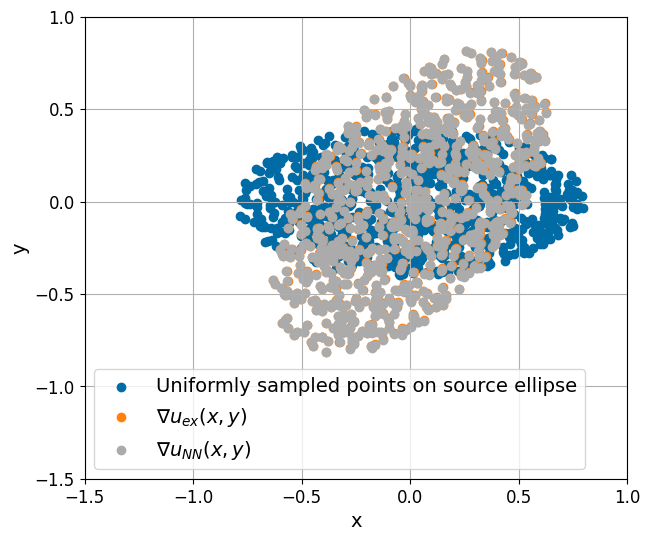} 
\end{center}
\caption{
Transport of a constant density from an ellipse domain to a rotated ellipse domain. 
Visualization of uniformly sampled source points and their image through exact and approximate transport maps.}
\label{fig:otsecond}
\end{figure}

\begin{figure}[ht!]
\begin{center}    
\includegraphics[width=0.5\linewidth]{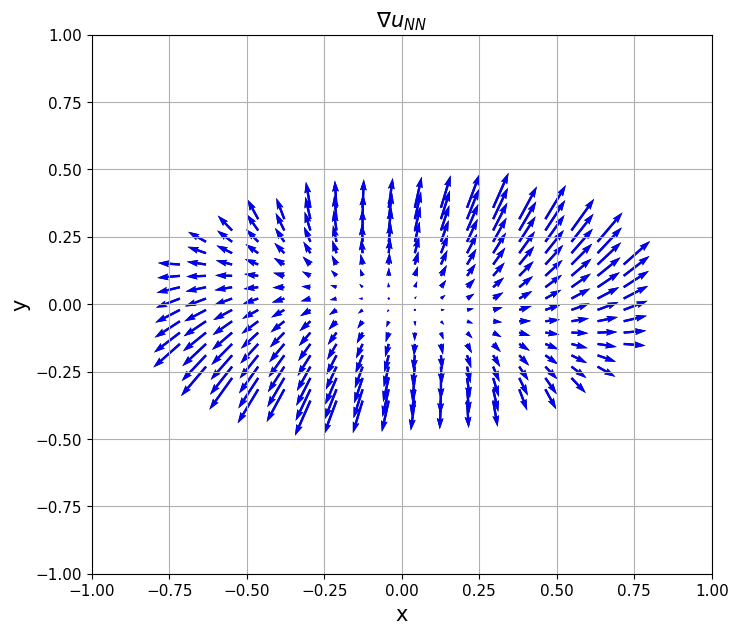} 
\end{center}
\caption{
Transport of a constant density from an ellipse domain to a rotated ellipse domain.
Visualization of the transport map $\nabla u_{NN}$.
}
\label{fig:otsecond_field}
\end{figure}

\begin{figure}[ht!]
\begin{center}    
\includegraphics[width=0.8\linewidth]{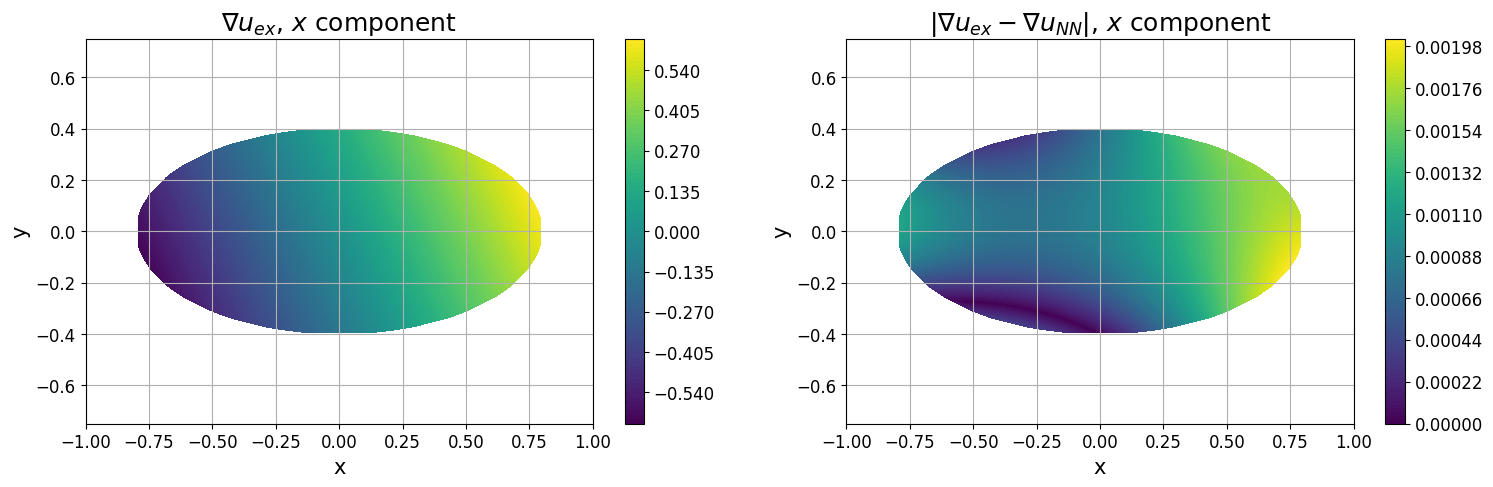} 
\includegraphics[width=0.8\linewidth]{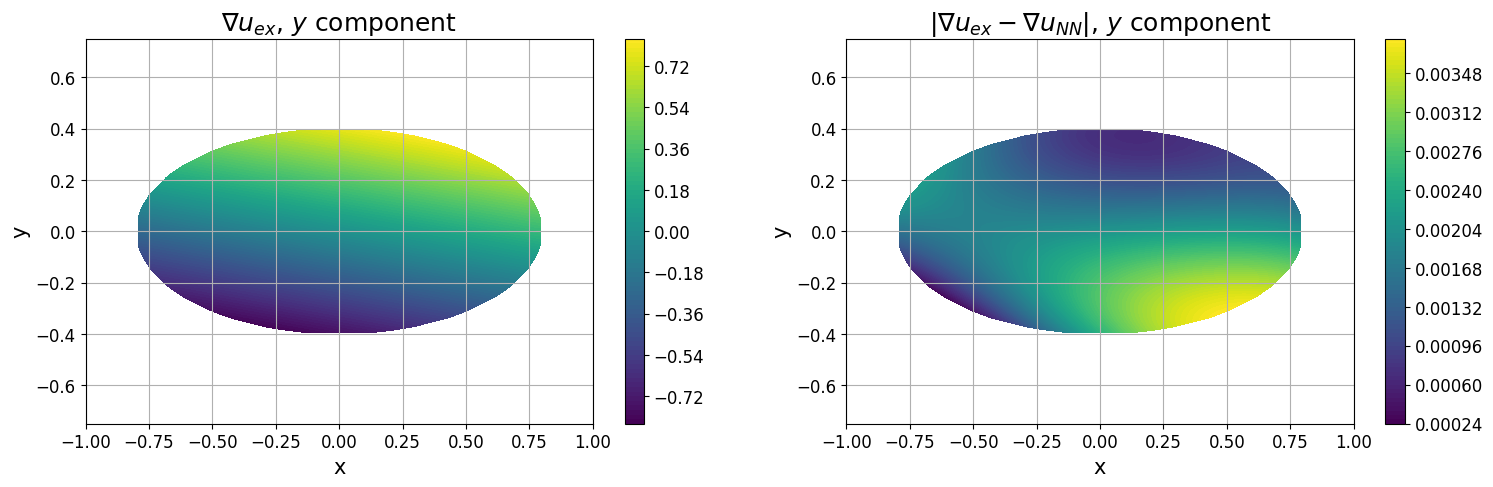} 
\end{center}    
    \caption{Transport of a constant density from an ellipse domain to a rotated ellipse domain. 
    Top left: first $x$-coordinate of the transport map $\nabla u_{NN}$; 
    Top right: Pointwise absolute error $\left| \nabla u_{NN,x} - \nabla u_{ex,x} \right|$;  
    Bottom left: second $y$-coordinate of the transport map $\nabla u_{NN}$; 
    Top right: Pointwise absolute error $\left| \nabla u_{NN,y} - \nabla u_{ex,y} \right|$. }
    \label{fig:otsecond_xerr_yerr}
\end{figure}

\begin{figure}[ht!]
\begin{center}    
\includegraphics[width=0.50\linewidth]{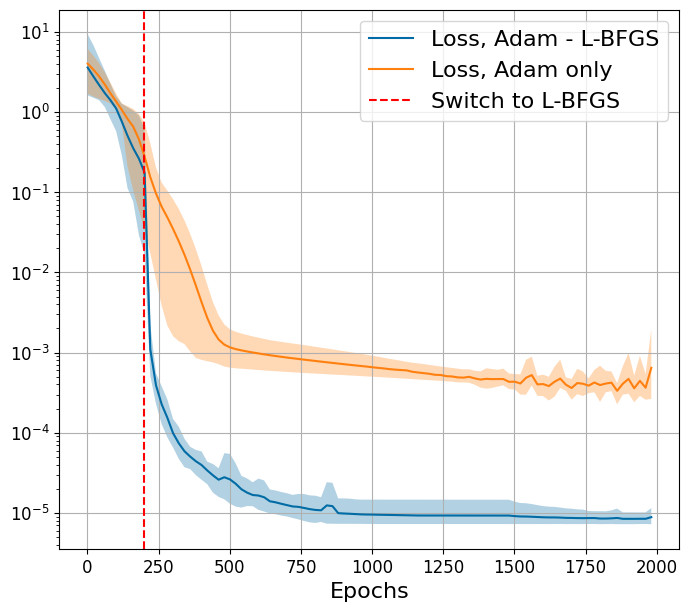} 
\end{center}    
    \caption{
    Transport of a constant density from an ellipse domain to a rotated ellipse domain. 
    Loss function vs epochs. 
    For the L-BFGS optimizer, the total epoch count is taken as the number of (outer) L‑BFGS iterations times the number of inner sub‐iterations. 
    Shaded bands show the 5th–95th percentile range over ten independent runs.
    }
    \label{fig:otsecond_loss}
\end{figure}

The second example aims at transporting a constant probability density function supported in an ellipse domain $\calX \subset \R^2$ to another constant probability density function whose support is a rotated ellipse domain $\calY \subset \R^2$.
More precisely, let $\mathcal{B}_1$ be the unit disk in $\R^2$ and $\mathcal{X} = M_{\mathcal{X}}\mathcal{B}_1$ and $\mathcal{Y} = M_{\mathcal{Y}}\mathcal{B}_1$, $M_{\mathcal{X}}, M_{\mathcal{Y}}$ being positive definite matrices. 
The optimal map can be explicitly calculated in $\mathbb{R}^2$ and is given by
\begin{equation*}
T({\bx})=M_{\mathcal{Y}}\mathcal{R}_{\theta}M_{\mathcal{X}}^{-1}{\bx}, \qquad {\bx}\in \mathcal{X},
\end{equation*}
where $\mathcal{R}_{\theta}$ is the rotation matrix
\begin{equation*}
\mathcal{R}_{\theta} = \begin{pmatrix}
\cos(\theta) & -\sin(\theta) \\
\sin(\theta)  & \cos(\theta) 
\end{pmatrix},
\end{equation*}
and the angle $\theta$ is calculated as follows:
\begin{equation*}
\tan(\theta) = \mathrm{trace}(M_{\mathcal{X}}^{-1}M_{\mathcal{Y}}^{-1}J)/\mathrm{trace}(M_{\mathcal{X}}^{-1}M_{\mathcal{Y}}^{-1}),
\end{equation*}
where 
\begin{equation*}
J = \begin{pmatrix}
0 & -1 \\
1  & 0
\end{pmatrix}.
\end{equation*}
We test the following example where
\begin{equation*}
M_{\mathcal{X}}= \begin{pmatrix}
0.8 & 0 \\
0  & 0.4
\end{pmatrix}
, \qquad
M_{\mathcal{Y}} = \begin{pmatrix}
0.8 & 0.2 \\
0.2  & 0.6
\end{pmatrix}.
\end{equation*}
We run the optimising process with $400$ Adam epochs and $1000$ L-BFGS epochs. 
Figure~\ref{fig:otsecond} shows how uniformly sampled points in the original ellipse $\calX$ are transported to the rotated ellipse $\calY$; it illustrates an appropriate, uniform, repartition of the sampled points and an accurate transport of those points, whose final position coincides well with their exact values. 
Figure~\ref{fig:otsecond_field} visualizes the corresponding vector field $\nabla u_{NN}$.
Figure~\ref{fig:otsecond_xerr_yerr} visualizes the approximation error on both components of $\nabla u_{NN}$, and highlights that the error is of the order of $10^{-3}$ pointwise. 
Finally, as in the previous test case, Figure~\ref{fig:otsecond_loss} illustrates the loss profiles. Again, the L‑BFGS optimizer converges more quickly and achieves lower loss values than the Adam optimizer.

\clearpage

\subsection{
Transport of a Gaussian density in the unit square to a uniform density in the unit square
}\label{sect:GtoU}
\begin{figure}[ht!]
\begin{center}
\includegraphics[width=\linewidth]{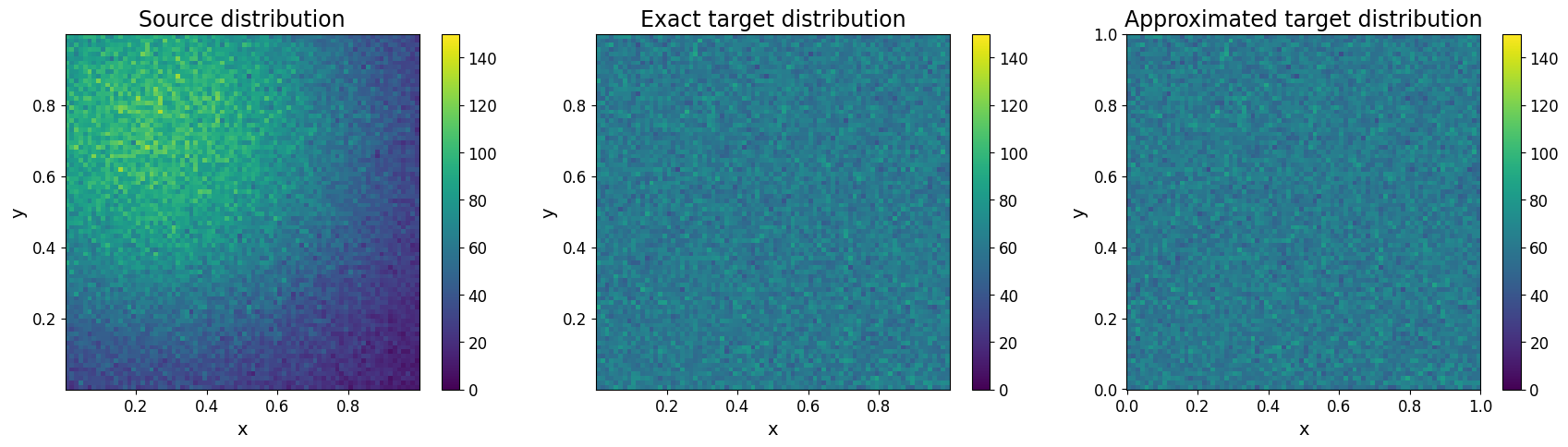} 
\end{center}
\caption{
Transport of a Gaussian density to a uniform density. 
Left: colormap obtained by the sampling of points according to the original Gaussian density (labeled here as source distribution); 
Middle: colormap obtained by the transport of the original density with the exact solution $\nabla u_{ex}$ (labeled here as exact target distribution); 
Right: colormap obtained by the transport of the original density with the approximated solution $\nabla u_{NN}$ (labeled here as approximated target distribution). 
($N = 100000$ points sampled into the cells of a structured grid of size $80 \times 80$).
}
\label{fig:otthird_2}
\end{figure}

In order to mimic the optimal transport of piles of debris \cite{Monge1781}, we transport the unit square domain with an initial Gaussian density function, namely 
\begin{equation}\label{eq:gaussian1}
f({\bx}) = c_0\exp\left\{-\dfrac{1}{2\sigma^2}({\bx}-{\bx}_0)^2\right\}\chi_{[0,1]^2}({\bx}), 
\end{equation}
where $c_0$ is a normalization constant, $\sigma^2 = 0.25$ and ${\bx}_0 = (0.25, 0.75)$. 
This function is transported into a uniform density on the same unit square domain:
$$
g({\bx}) = \chi_{[0,1]^2}({\bx}).
$$
Figure~\ref{fig:otthird_2} illustrates the approximated optimal transport map. 
The histogram in Figure~\ref{fig:otthird_2}~(left) shows the distribution of points sampled from the Gaussian distribution (\ref{eq:gaussian1}) in the unit square.
Then the points are transported to a uniform distribution through the exact map $\nabla u$ (center) and the approximated map $\nabla u_{NN}$ (right).

\begin{figure}[ht!]
\begin{center}    
\includegraphics[width=0.75\linewidth]{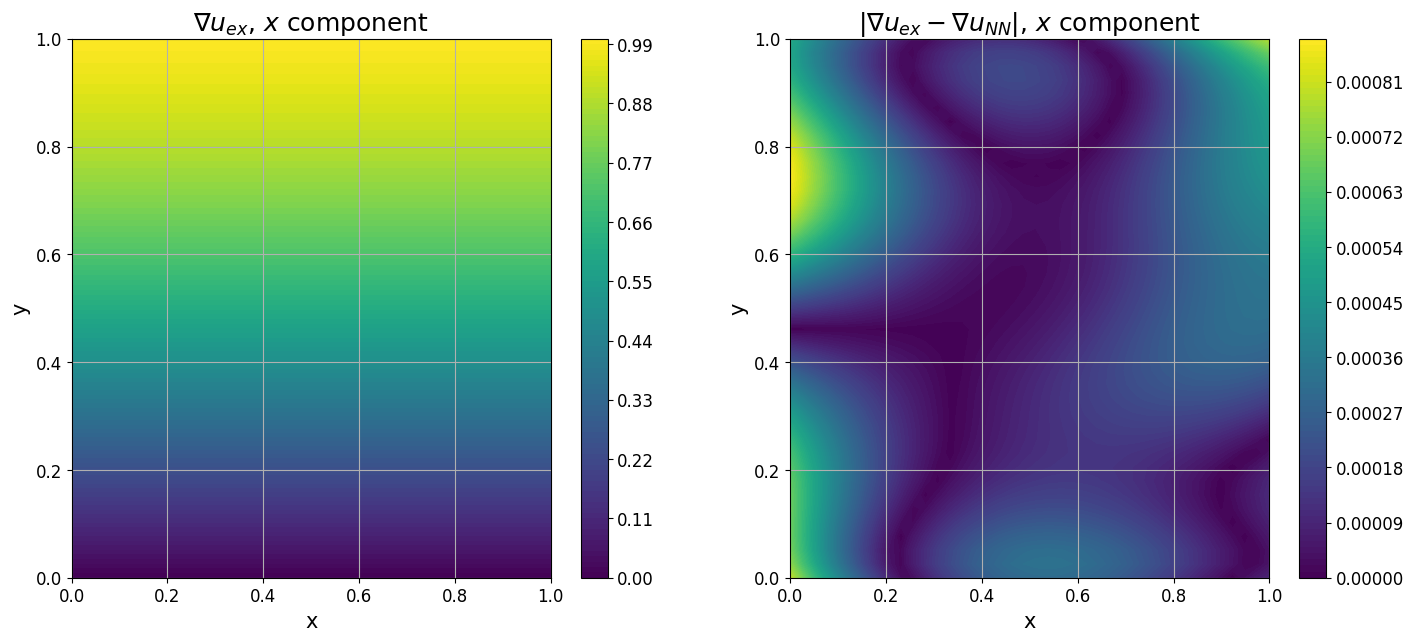} 
\includegraphics[width=0.75\linewidth]{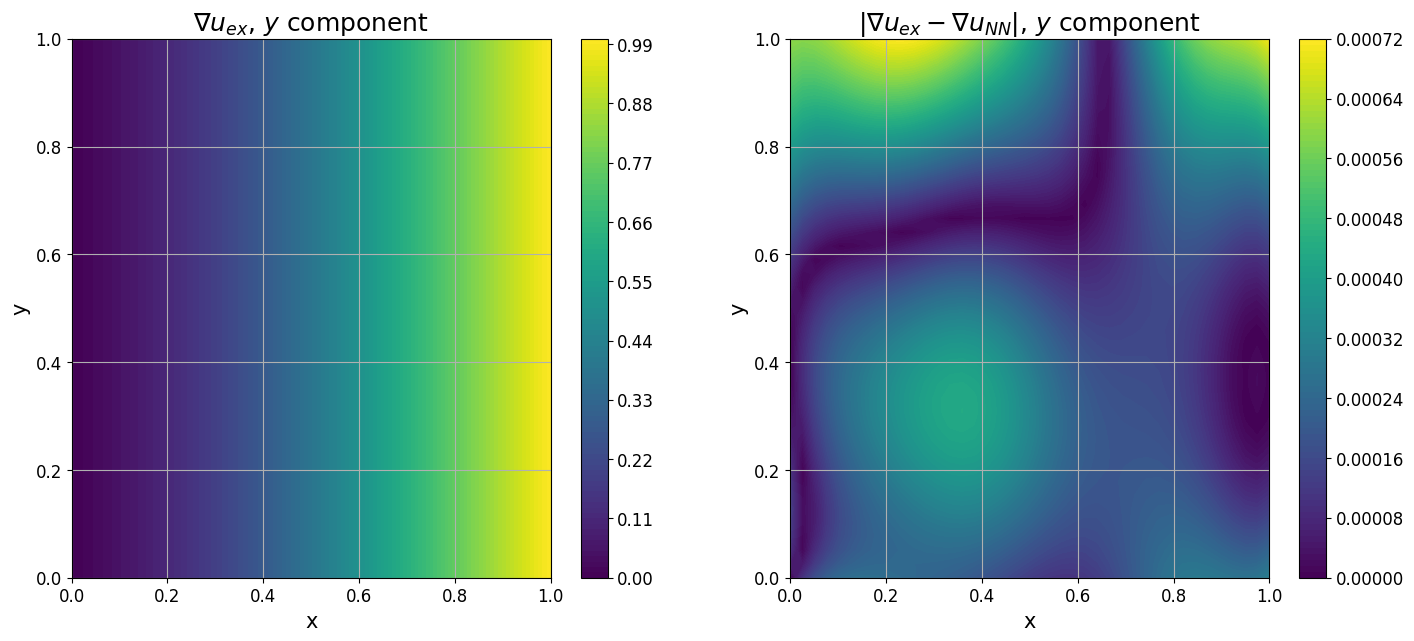} 
\end{center}    
    \caption{
    Transport of a Gaussian density to a uniform density.  
    Top left: first $x$-coordinate of the transport map $\nabla u_{NN}$; 
    Top right: Pointwise absolute error $\left| \nabla u_{NN,x} - \nabla u_{ex,x} \right|$;  
    Bottom left: second $y$-coordinate of the transport map $\nabla u_{NN}$; 
    Top right: Pointwise absolute error $\left| \nabla u_{NN,y} - \nabla u_{ex,y} \right|$. }
    \label{fig:otthird_xerr_yerr}
\end{figure}

We can observe that the transported solution re-dispatches appropriately uniformly the points on the grid, and that the exact and the approximated solutions are very similar to each other.
Figure~\ref{fig:otthird_xerr_yerr} visualizes the approximation error on both components of $\nabla u_{NN}$. The error is still of a similar order, and we can observe that it concentrates in the corners of the unit square domain. 

Note that, to illustrate the original and approximated solutions in Figure~\ref{fig:otthird_2} and in the remainder of the experiments, we choose to randomly generate a sample of visualization points following the original probability distribution $f$ on the domain $\mathcal{X}$. We then illustrate those points at their original location to visualize the original density function, together with the image of those points via the (approximated) transport map to obtain the final points. We then illustrate those final points to visualize the final density function.

\subsection{
Transport of a Gaussian density in the unit square to another Gaussian density in the unit square 
}
\begin{figure}[ht!]
\begin{center}
\includegraphics[width=0.95\linewidth]{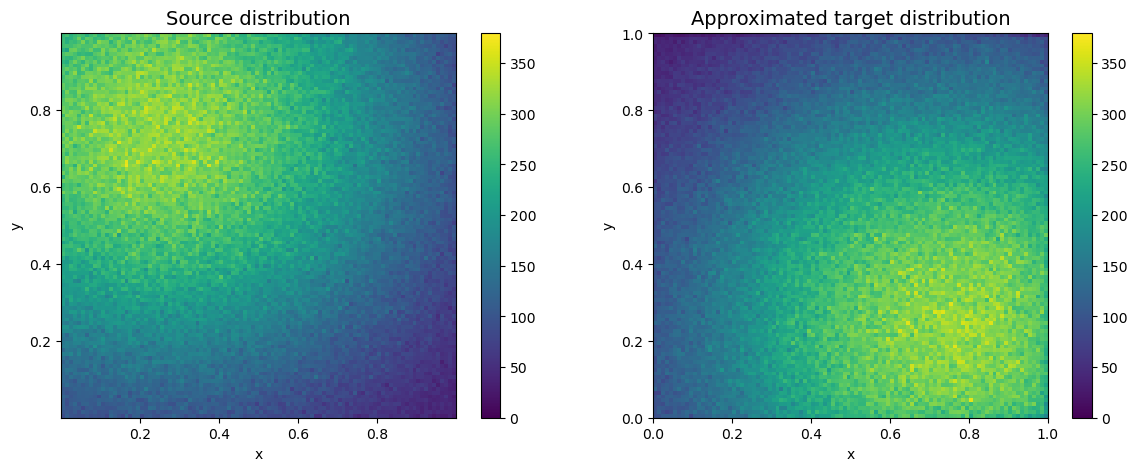} 
\end{center}
\caption{
Transport of a Gaussian density to another Gaussian density. 
Left: colormap obtained by the sampling of points according to the original Gaussian density (labeled here as source distribution); 
Right: colormap obtained by the transport of the original density with the approximated solution $\nabla u_{NN}$ (labeled here as approximated target distribution). 
($N = 100000$ points sampled into the cells of a structured grid of size $80 \times 80$).
}
\label{fig:otgaussian_to_gaussian}
\end{figure}

We consider the transport of an initial Gaussian density function within the unit square domain into another Gaussian density function, in order to evaluate the numerical diffusion of the method when reconstructing a non-constant target density function. 
More precisely, we consider: 
\begin{equation*}
f({\bx}) = c_0\exp\left\{-\dfrac{1}{2\sigma^2}({\bx}-{\bx}_0)^2\right\}\chi_{[0,1]^2}({\bx}). 
\end{equation*}
where $c_0$ is a normalization constant, $\sigma^2 = 0.25$ and ${\bx}_0 = (0.25, 0.75)$, and 
$g({\bx})$ defined similarly, with $\sigma^2 = 0.25$ and ${\bx}_0 = (0.75, 0.25)$. 
The histogram in Figure~\ref{fig:otgaussian_to_gaussian}~(left) shows the distribution of points sampled from the source Gaussian distribution (\ref{eq:gaussian1}) in the unit square and then transported to the target Gaussian distribution through the map $\nabla u_{NN}$ (right). 
Even if there is no analytical solution, we observe that the distribution of the transported points very closely resembles a Gaussian distribution.

\subsection{
Transport of a bimodal Gaussian density in the unit square to a uniform density in the unit square
}

Finally, in an effort to mimic several suppliers transporting material from several sources, we present an example of multiple Gaussian sources being scattered into a uniform density field. 
We consider an anisotropic bimodal Gaussian distribution on the unit square, and a target uniform distribution on the same unit square, namely we define, with $\bx=(x_1,x_2)$: 

\begin{eqnarray}
f({\bx}) &=& c_0\left(\exp\left\{-\dfrac{1}{2\sigma^2_{xx}}(x_1-0.5)^2-\dfrac{1}{2\sigma^2_{yy}}(x_2-0.2)^2\right\} \right. 
\nonumber \\ && 
\left. + \exp\left\{-\dfrac{1}{2\sigma^2_{xx}}(x_1-0.5)^2-\dfrac{1}{2\sigma^2_{yy}}(x_2-0.8)^2\right\}\right)
\chi_{[0,1]^2}({\bx}), \label{eq:gaussian2}
\end{eqnarray}
where $c_0$ is a normalization constant, $\sigma^2_{xx} = 0.25$ and $\sigma^2_{yy} = 0.015625$, and
\begin{eqnarray*}
g({\bx}) = \chi_{[0,1]^2}({\bx}).
\end{eqnarray*}

\begin{figure}[ht!]
\begin{center}
\includegraphics[width=0.95\linewidth]{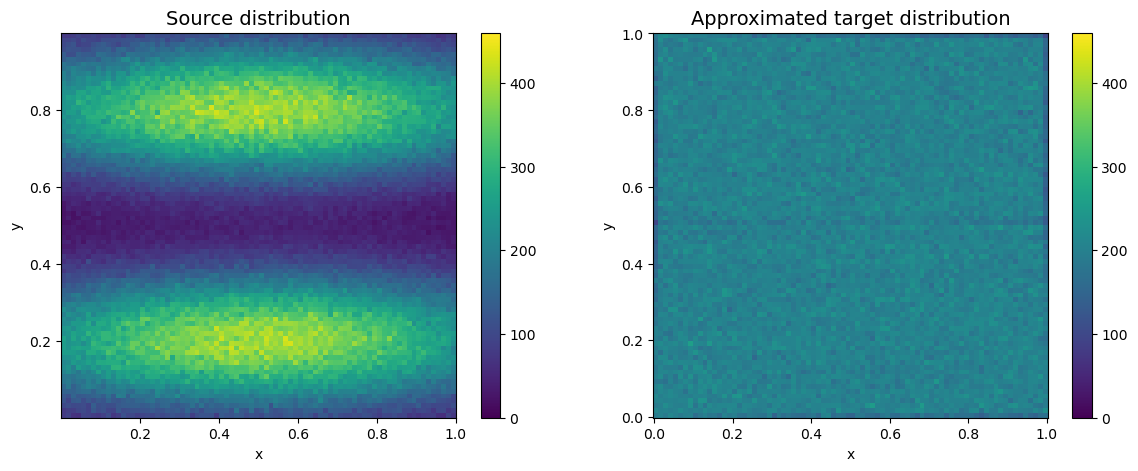} 
\end{center}
\caption{
Transport of a bimodal Gaussian density to a uniform density.
Left: colormap obtained by the sampling of points according to the original bimodal Gaussian density (\ref{eq:gaussian2}) (labeled here as source distribution); 
Right: colormap obtained by the transport of the original density with the approximated solution $\nabla u_{NN}$ (labeled here as approximated target distribution). 
($N = 100000$ points sampled into the cells of a structured grid of size $80 \times 80$).
}
\label{fig:otfourth_2}
\end{figure}

The histogram in Figure~\ref{fig:otfourth_2}~(left) shows the distribution of points sampled from the source bimodal Gaussian distribution (\ref{eq:gaussian2}) in the unit square and then transported to the uniform distribution through the map $\nabla u_{NN}$ (right). 
Even if the exact solution $u$ is unknown, we observe that the distribution of the transported points is uniform.

\subsection{
Transport of a 3D Gaussian density in the unit cube to a uniform density in the unit square
}
\begin{figure}[ht!]
\begin{center}
\includegraphics[width=0.95\linewidth]{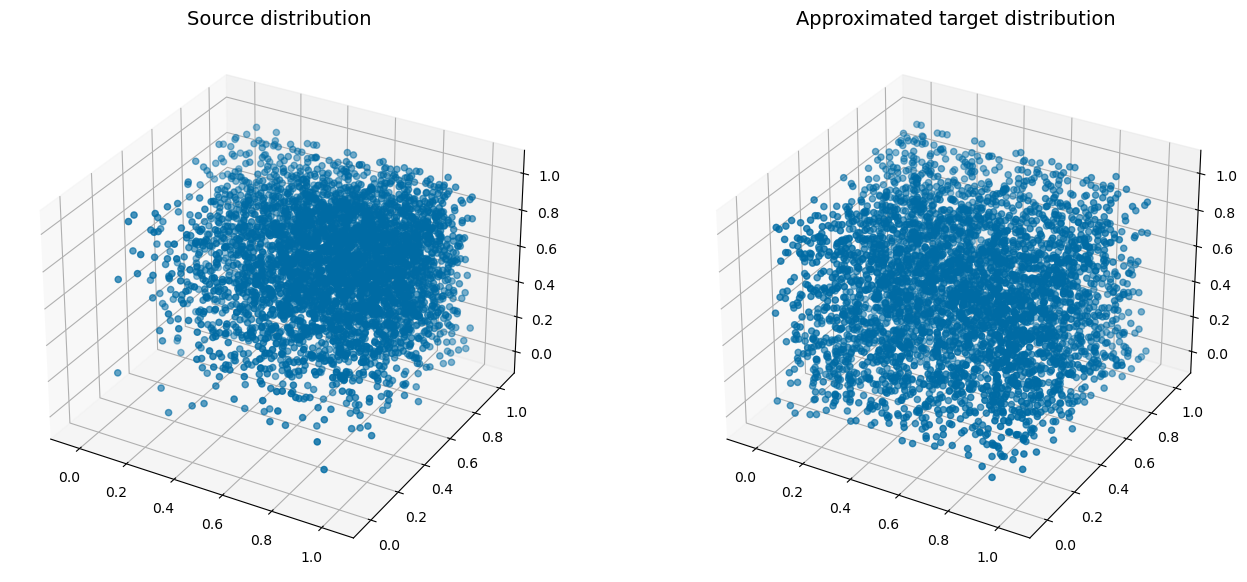} 
\end{center}
\caption{
Transport of a 3D Gaussian density to a uniform density.
Left: sampling of points inside the unit cube $\calX$ according to Gaussian density (\ref{eq:gaussian1_3d});
Right: location of the inside points transported with the approximated solution $\nabla u_{NN}$. 
}
\label{fig:otfourth_2_3D}
\end{figure}

\begin{figure}[ht!]
\begin{center}
\includegraphics[width=0.95\linewidth]{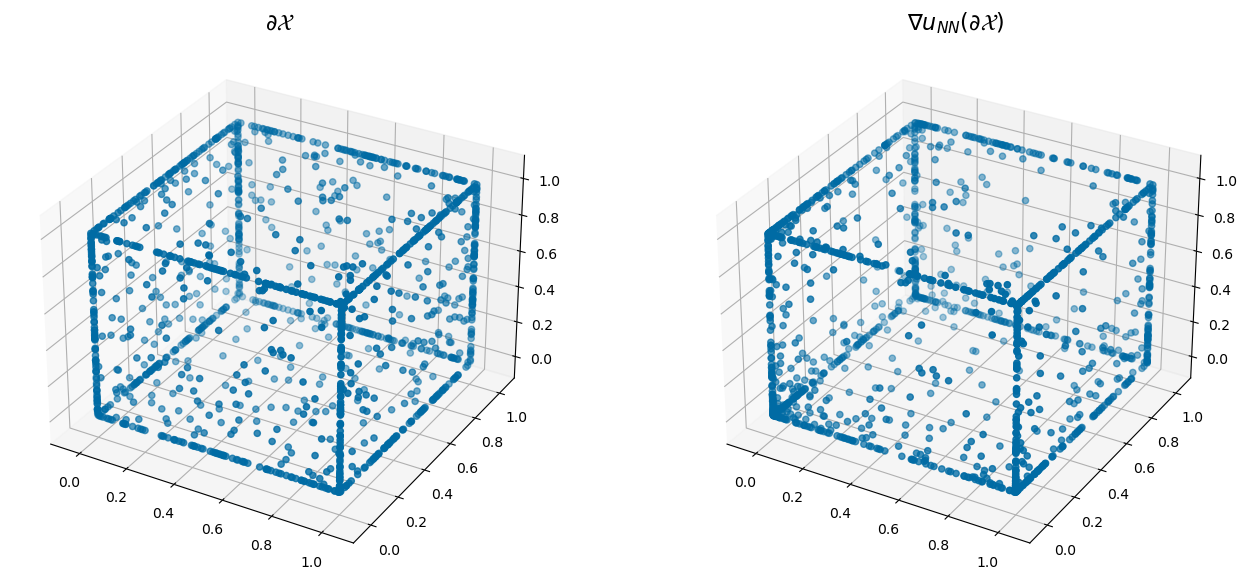} 
\end{center}
\caption{
Transport of a 3D Gaussian density to a uniform density.
Left: sample of points on the boundary of the unit cube $\partial \calX$, uniformly sampled on the boundary
Right: location of the boundary points transported with the approximated solution $\nabla u_{NN}$. 
}
\label{fig:otfourth_3}
\end{figure}

Finally, to illustrate the ability of the algorithm to be effective in more than two space dimensions, the test case in Section~\ref{sect:GtoU} is extended to the 3D case. 
We consider 

\begin{equation}\label{eq:gaussian1_3d}
f({\bx}) = c_0\exp\left\{-\dfrac{1}{2\sigma^2}({\bx}-{\bx}_0)^2\right\}\chi_{[0,1]^3}({\bx}), 
\end{equation}
where $c_0$ is a normalization constant, $\sigma^2 = 0.125$ and ${\bx}_0 = (0.75, 0.75, 0.75)$. 
This function is transported into a uniform density on the unit cube:
$$
g({\bx}) = \chi_{[0,1]^3}({\bx}). 
$$
Figure~\ref{fig:otfourth_2_3D} illustrates the approximated optimal transport map. 
We can observe that the transported solution re-dispatches appropriately uniformly the points on the grid. 
Figure~\ref{fig:otfourth_3} illustrates the transport of uniformly randomly sampled boundary points. 
We observe that $\nabla u_{NN}(\partial \calX) = \partial \calY$ as well. 

\section{Conclusions and Perspectives}\label{sect:conclusion}

A model based on a generalized Monge-Amp\`ere equation with transport boundary conditions has been considered to solve the optimal transport problem. 
A convex input PINNs method has been advocated for the numerical approximation of its convex solution. 
The loss function has been modified to account for the transport boundary conditions by using a discrete version of the Hausdorff distance. 

Since we look for a convex solution of the Monge-Amp\`ere equation, input convex neural networks are suitable to obtain appropriate approximations as they are seeking for solutions in the right subspaces.  
This advantage comes at the cost of increased training complexity. 

A sensitivity analysis has shown that the model is not particularly affected by the random sampling of the collocation and boundary points. 
Numerical results have focused on transporting uniform and Gaussian distribution functions in two dimensions of space. 
However, we have shown that results can be extended to higher dimensions. 

The relatively high computational expense currently limits the practical problem size and network scale. Furthermore the variability of the results over the samples show that the error coming from the PINNs approach is rather difficult to control.
Thus, future avenues of investigations should therefore focus on reducing the CPU overhead of the method, and controlling the error of the method, for instance with hybrid methods. 
Additional perspectives also include the further incorporation of several sources and targets, via multimodal Gaussian density functions, for real-life applications such as raw material transportation with several sources, as well as the treatment of obstacles. 
Finally, the investigation of the Monge-Amp\`ere in larger number of space dimensions could be addressed. 

\section*{Acknowledgements}
The authors thank Archie Lawson (EPFL) for implementation support during an EPFL Bachelor semester project, and Prof. Marco Picasso (EPFL) for fruitful discussions.



\end{document}